\documentclass[12pt]{amsart}
\usepackage[all]{xy}
\usepackage{relsize}
\usepackage{amssymb}
\usepackage{amsthm}
\usepackage{hyperref}
\usepackage{anyfontsize}
\hypersetup{colorlinks=true,linkcolor=blue,citecolor=magenta}
\usepackage{amsmath}
\usepackage{amscd,enumitem}
\usepackage{verbatim}
\usepackage{eurosym}
\usepackage{float}
\usepackage{color}
\usepackage{dcolumn}
\usepackage[mathscr]{eucal}
\usepackage[all]{xy}
\usepackage{hyperref}
\usepackage{bbm}
\usepackage[textheight=8.8in, textwidth=6.8in]{geometry}
\usepackage{mathtools}

\newtheorem*{conj*}{Conjecture}

\newtheorem{theorem}{Theorem}[section]
\newtheorem{remark}[theorem]{Remark}

\newtheorem{prop}[theorem]{Proposition}

\newtheorem*{notation*}{Notation}

\newtheorem{remarks}[theorem]{Remarks}

\newtheorem{definition}[theorem]{Definition}

\newtheorem{example}{Example}[section]

\newcommand{\Z}{\mathbb{Z}}

\newcommand{\R}{\mathbb{R}}
\newcommand{\SL}{\operatorname{SL}}
\newcommand{\C}{\mathbb{C}}

\renewcommand{\H}{\mathbb{H}}

\renewcommand{\b}{\beta}
\newcommand{\g}{\gamma}
\newcommand{\G}{\Gamma}
\renewcommand{\t}{\tau}
\renewcommand{\r}{\rho}
\renewcommand{\l}{\lambda}

\renewcommand{\arg}{\mathrm{Arg}}

\newcommand{\cro}{C_\rho}

\newcommand{\F}[1]{{}_2 F_1\left( #1 \right)}

\newcommand{\ol}[1]{\overline{#1}}

\numberwithin{equation}{section}

\usepackage{tikz}

\title{
Inversion Formulas for the $j$-function around Elliptic Points
} 
\author{Alejandro De Las Penas Castano}
\address{Department of Mathematics, University of Virginia, Charlottesville, VA 22904}
\email{ad7ag@virginia.edu}
\author{Badri Vishal Pandey}
\address{Department of Mathematics, University of Virginia, Charlottesville, VA 22904}
\email{bp3aq@virginia.edu}

\begin{document}

\begin{abstract} 
Recently, Hong, Mertens, Ono, and Zhang \cite{hong2021proof} proved a conjecture of C\u{a}ld\u{a}raru, He, and Huang \cite{caldararu2021moonshine} that expresses the Taylor series of the modular $j$-function around the elliptic points $i$ and $\r=e^{\pi i/3}$ as rational functions arising from the signature 2 and 3 cases of Ramanujan's theory of elliptic functions to alternative bases. We extend these results and give inversion formulas for the $j$-function around $i$ and $\r$ arising from Gauss' hypergeometric functions and Ramanujan's theory in signatures 4 and 6.
\end{abstract}

\maketitle
\section{Introduction and statement of results}

The Klein $j$-function
\[
    j(\t):=\frac{1}{q}+744+196884q+21493760q^2+\cdots\qquad(q=e^{2\pi i\t},\; \t\in\H)
\]
is a modular function on the full modular group $\SL_2(\Z)$. It is of great importance to number theory. In the theory of elliptic curves, the $j$-function parametrizes isomorphism classes of elliptic curves over $\C$. In Class Field Theory, its values at CM points, the so called singular moduli, generate Hilbert Class Fields of imaginary quadratic extensions. Another famous example of its importance is the observation that the Fourier coefficients of the $j$-function encode the graded dimensions of the infinite dimensional graded algebra of the Monster group. This observation led to the Monstrous Moonshine conjecture and its eventual proof by Borcherds \cite{Borcherds92}.

The $j$-function defines a bijective holomorphic function from $\SL_2(\Z)\backslash\H$ to $\C$. In particular, the $j$-function has an inverse function. Due to its central role in the theories stated above, among many others, it is natural to seek explicit formulas for its inverse map. Work in this direction began with Ramanujan giving striking formulas for $1/\pi$ which gave rise to Ramanujan's theory of elliptic functions to alternative basis (cf. \cite{Ramanujan1914}, \cite{Berndt1995}, and \cite{berndt_chan_1999}). This theory produced several explicit formulas which express $j(\t)$ as a rational function in $t$, where the parameters $\t$ and $t$ are related via Gauss' hypergeometric functions (see Section  \ref{sec:hypergeometric}). 

More recently, a conjecture of C\u{a}ld\u{a}raru, He, and Huang \cite{caldararu2021moonshine} cast new light on the problem of inverting the $j$-function. In contrast to Ramanujan's Theory which uses the Fourier expansion of the $j$-function around the cusp $i\infty$ \cite{berndt_chan_1999}, their conjecture is about the Taylor series expansion around the elliptic points $i$ and $\rho:=e^{\pi i/3}$. They conjecture that these Taylor series, when specialized at the normalized flat coordinate of the corresponding moduli space of versal deformations of elliptic curves, are the rational functions that appear in the classical hypergeometric inversion formulae for the $j$-function. Shortly after, Hong, Mertens, Ono, and Zhang \cite{hong2021proof} proved their conjecture.

To state their results more precisely, let $\H$ be upper half plane, and $\mathbb{D}$ be the unit disc. For $\t_*\in\{i,\r\},$ we have the uniformizing map $S_{\t_*}:\H\to\mathbb{D}$ and its inverse $S_{\t_*}^{-1}:\mathbb{D}\to\H$ defined by 
$$S_{\t_*}(\t)=\frac{\t-\t_*}{\t-\ol{\t_*}} \quad\text{ and }\quad S^{-1}_{\t_*}(w)=\frac{\t_*-\ol{\t_*}w}{1-w}.$$
Next, we renormalize the conformal map
  \[
    s^{-1}_{\t_*}(w):=S_{\t_*}^{-1}\left(\frac{w}{2\pi\Omega_{\r}^2}\right),
  \]
where $\Omega_{\t_*}$ is the standard Chowla-Selberg periods (for example, see (96) of Section 6.3 of \cite{Zagier2008}) defined by
\[
    \Omega_{i}:=
                    \dfrac{1}{\sqrt{8\pi}}\dfrac{\G(1/4)}{\G(3/4)}\quad\text{and}\quad  \Omega_{\r}:=
                    \dfrac{1}{\sqrt{6\pi}}\left(\dfrac{\G(1/3)}{\G(2/3)}\right)^{3/2}.
\]
With these uniformizing maps, the Taylor series of the $j$-function around $i$ and $\r$ are defined as:
\begin{align*}
    j(s_i^{-1}(w))&=1728+20736w^2+105984w^4+\frac{1594112}{5}w^6+\cdots\\
    j(s_\r^{-1}(w))&=13824w^3-39744w^6+\frac{1920024}{35}w^9-\frac{1736613}{35}w^{12}+\cdots.
\end{align*}
These formulas follow from the theory of Taylor coefficients of modular forms (see for example Section 5.4 of \cite{cohen2017modular}).

Hong, Mertens, Ono, and Zhang considered the two distinguished power series $c_i(t)$ and $c_\r(t)$ of C\u{a}ld\u{a}raru, He, and Huang (see Section 2.1 of \cite{hong2021proof}) whose first few terms are
\begin{align*}
    c_i(t)&=t+t^3+\frac{32}{15}t^5+\frac{17}{3}t^7+\frac{1054}{63}t^{9}+\frac{368} {7}t^{11}+\frac{4652300}{27027}t^{13}+\cdots &(|t|<1/2)\\
    c_\r(t)&=t-\frac{1}{3}t^4+\frac{103}{315}t^7-\frac{169}{405}t^{10}+\frac{522169}{868725}t^{13}-\frac{186119}{200475}t^{16}+\cdots  &(|t|<1).
\end{align*}
They proved that surprisingly the Taylor series of the $j$-function around $i$ and $\r$ evaluated at $c_i(t)$ and $c_{\r}(t)$ respectively, turn out to be rational functions in $t$. Namely, they proved that
\[
    j(s_i^{-1}(c_i(t))=64\frac{(3+4t)^3}{(1-4t^2)^2}\quad
\text{ and }\quad
    j\left(\frac{s_\r^{-1}(c_\r(t))+1}{3}\right)=27t^3\left(\frac{8-t^3}{1+t^3}\right)^3.
\]
Their proof used the theory of hypergeometric functions of signatures 2 and 3. Namely, they realized $c_i$ and $c_\rho$ as quotients of hypergeoemetric functions. 

In view of this, it is natural to ask whether there are other examples of this phenomenon. More precisely, does the theory of hypergeometric functions in signature 4 and 6 yield other inversion formulas for the $j$-function around the elliptic points $i$ and $\r$? We answer this question here.

Let us define
\begin{equation}\label{eq:def-Ci-Cr}
    C_i(t):=t\cdot\frac{\F{\frac{3}{4},\frac{3}{4};\frac{3}{2};4t^2}}{\F{\frac{1}{4},\frac{1}{4};\frac{1}{2};4t^2}}\quad\text{and}\quad \cro(t):=\dfrac{t^2}{2}\dfrac{\F{\tfrac{5}{6},\tfrac{5}{6};\tfrac{5}{3};-2t^3}}{\F{\tfrac{1}{6},\tfrac{1}{6};\tfrac{1}{3};-2t^3}},
\end{equation}
where ${}_{2}F_1$ is Gauss' hypergeometric function (see Section \ref{sec:hypergeometric}). Then we have the following theorems:

\begin{theorem}\label{thm:sig4}
  If $|t|<1/2$, then we have
  \[
    j\left(\frac{s^{-1}_i(C_i(t))+1}{2}\right)=64\frac{(16t^2-3)^3}{4t^2-1}.
  \]
\end{theorem}

\begin{theorem}\label{thm:sig6}
  If $|t|<1/\sqrt[3]{2}$, then we have
  \[
    j\left(s^{-1}_{\r}(\cro(t))\right)=-\dfrac{1728t^6}{2t^3+1}.
  \]
\end{theorem}

\begin{remarks}$\;$
\begin{enumerate}
    \item The above formulas can be thought of as inversion formulas of the $j$-function around $i$ and $\r$ because when specialized to $t=0$, the above formulas reduce to the classic identities
\[
    j(i)=1728\quad\text{and}\quad j(\r)=0.
\]

    \item By analytic continuation, we can extend the domain where the formulas in Theorems \ref{thm:sig4} and \ref{thm:sig6} are valid to wherever $C_i(t)$ and $\cro(t)$ are defined, see the comments at the end of section \ref{sec:pf-thm1.2}. However, for explicit computations, we require $|t|<1/2$ and $|t|<1/\sqrt[3]{2}$ in order to compute the power series expansions of $C_i$ and $\cro$ respectively.
    
    \item These results tell us that finding an approximate solution to $j(\t)=\alpha$ boils down to solving a degree six polynomial equation, see Subsection \ref{sec:examples} for examples.
\end{enumerate}
\end{remarks}

In Section \ref{sec:nuts-and-bolts} we review certain identities involving the hypergeometric functions ${}_2F_1$ and the inversion formulas mentioned above. In Sections \ref{sec:pf-thm1.1} and \ref{sec:pf-thm1.2}, we prove Theorems \ref{thm:sig4} and \ref{thm:sig6} respectively. Finally, in Section \ref{sec:examples} we offer some explicit examples of these results.

\section*{Acknowledgement}
We wish to thank Ken Ono for helpful discussions related to this paper, and for providing research support with his NSF Grant DMS-2055118.

\section{Nuts and bolts}\label{sec:nuts-and-bolts}

\subsection{Hypergeometric Functions}\label{sec:hypergeometric}Here we recall the necessary facts about hypergeometric functions.

\begin{definition}
  Let $a,b\in\R$ and $c\in \R\setminus\Z^-$. The Gaussian hypergeometric function is defined as:
  \[
    \F{a,b;c;z}:=\sum_{n=0}^\infty \frac{(a)_n(b)_n}{(c)_n n!}z^n,
  \]
  for $|z|<1$, where $(s)_n$ is the Pochhammer symbol defined as
  \[
    (s)_n=s(s+1)\cdots(s+n-1).
  \]
  \end{definition}
  \begin{notation*}
    For a positive integer $r$, we define the following special hypergeometric function
  \begin{equation}\label{eq:lambda}
    \l_r(z):=\F{\frac{1}{r},1-\frac{1}{r};1,z}.
  \end{equation}
  \end{notation*}

\begin{remarks}$\,$
\begin{enumerate}\label{rm:discont}
  \item The above series converges absolutely and uniformly on compact sets in the unit disk $|z|<1$ and using its integral representation it can be extended analytically to the region $\C\setminus[1,\infty)$, i.e., for $|\arg(1-z)|<\pi$ (see for example \cite{Erdelyi1955}). On the line $[1,\infty)$, we extend the definition as:
  \[
    \F{a,b;c;x}:=\lim_{\varepsilon\rightarrow 0^+} \F{a,b;c;x+i\varepsilon}.
  \]
  Note that this makes ${}_2F_1$ well-defined on the whole complex plane, analytic on the region $|\arg(1-z)|<\pi$ but discontinuous on the line $[1,\infty)$. Furthermore, the difference between the principal branches on the two sides of the branch cut is:
  \begin{align*}
    &\lim_{\varepsilon\rightarrow 0^+} \F{a,b;c;x+i\varepsilon}-\lim_{\delta\rightarrow 0^+} \F{a,b;c;x-i\delta}\\
    &=\frac{2\pi i\G(c)}{\G(a)\G(b)\G(c-a-b+1)}(1-x)^{c-a-b}\F{c-a,c-b;c-a-b+1;1-x}
  \end{align*}
  (see Section 15.2 of \cite{NIST:DLMF})
  \item   Note that in the definition of ${}_2F_1$, the arguments $a$ and $b$ are symmetric so we will swap them as the situation requires without additional comments.
\end{enumerate}
\end{remarks}

We require the following classical hypergeometric transformation law.

\begin{prop}\label{prop:linear-cov}
  (Equation 15.10.33 of \cite{NIST:DLMF})For $0<|\arg(1-z)|<\pi$, we have
  \begin{multline*}
    \F{a,b;c;z}
    =\frac{\G(1-b)\G(c)}{\G(a-b+1)\G(c-a)}\left(\frac{1}{z}\right)^a \F{a-c+1,a;a-b+1;\frac{1}{z}}\\
    +\frac{\G(1-b)\G(c)}{\G(a)\G(c-a-b+1)}\left(1-\frac{1}{z}\right)^{c-a-b}\left(-\frac{1}{z}\right)^b \F{c-a,1-a;c-a-b+1;1-\frac{1}{z}}.
  \end{multline*}
\end{prop}

\begin{remark}\label{rm:identity-for-x>1}
  In view of Remark \ref{rm:discont}, the above proposition can be extended to $x>1$ in the following manner:
  \begin{multline*}
    \F{a,b;c;x}
    =\frac{\G(1-b)\G(c)}{\G(a-b+1)\G(c-a)}\left(\frac{1}{x}\right)^a \lim_{\delta\rightarrow 0^+}\F{a-c+1,a;a-b+1;\frac{1}{x}-i\delta}\\
    +\frac{\G(1-b)\G(c)}{\G(a)\G(c-a-b+1)}\left(1-\frac{1}{x}\right)^{c-a-b}\left(-\frac{1}{x}\right)^b \F{c-a,1-a;c-a-b+1;1-\frac{1}{x}}.
  \end{multline*}
\end{remark}

\subsection{Inversion Formulas} Here we recall two classical inversion formulas for $j$-function in terms of $_{2}F_1$ hypergeometric functions.

\begin{prop}\label{thm:j-inv-q4}(Theorem 9.5-6 of \cite{Berndt1995})
If $\t\in\H$ and $\g$ satisfies
  \[
    \t=\frac{i}{\sqrt{2}}\frac{\F{\tfrac{1}{4},\tfrac{3}{4};1;1-\g}}{\F{\tfrac{1}{4},\tfrac{3}{4};1;\g}},
  \]
  then we have
  \begin{align}
    j(\t)&=\frac{64(1+3\g)^3}{\g(\g-1)^2}\label{eq:inv-q4-j}.
  \end{align}
\end{prop}
\begin{remark}
  In \cite{Berndt1995}, $j(\t)$ is not explicitly calculated, but Theorems 9.5 and 9.6 of \cite{Berndt1995} give formulas for $E_4(\t)$ and $E_6(\t)$, respectively, from which the formula for $j(\t)$ is immediately deduced.
\end{remark}

\begin{prop}\label{thm:j-inv-q6}(Theorem 11.4-5 of \cite{Berndt1995})
  Let $\g\in\C$ and $\b=-\g^3/2$. If we define
  \[
    \t=i\cdot\frac{\F{\tfrac{1}{6},\tfrac{5}{6};1;1-\b}}{\F{\tfrac{1}{6},\tfrac{5}{6};1,\b}},
  \]
  then we have
  \begin{align}
    j(\t)&=\frac{1728}{1-(1-2\b)^2}=\frac{-1728}{\g^3(2+\g^3)}\label{eq:inv-q6-j}.
  \end{align}
\end{prop}

\section{Proof of Theorem \ref{thm:sig4}}\label{sec:pf-thm1.1}
The following calculations depend on the argument of $t$. For the moment, we assume that $0<\arg(t)<\tfrac{\pi}{2}$. Apply Proposition \ref{prop:linear-cov} to $a=b=\tfrac{3}{4}$, $c=\tfrac{3}{2}$ and $z=4t^2$ to $\F{\tfrac{3}{4},\tfrac{3}{4};\tfrac{3}{2};4t^2}$, using the functional equation $s\G(s)=\G(s+1)$, and the lambda notation in \eqref{eq:lambda} we get
  
  \begin{align}
    \F{\frac{3}{4},\frac{3}{4};\frac{3}{2};4t^2}
    &=\frac{\G\left(\tfrac{1}{4}\right)\G\left(\tfrac{3}{2}\right)}{\G\left(\tfrac{3}{4}\right)}\left(\frac{1}{4t^2}\right)^{3/4}\left[ \F{\frac{1}{4},\frac{3}{4};1;\frac{1}{4t^2}}+e^{3\pi i/4}\F{\frac{1}{4},\frac{3}{4};1;1-\frac{1}{4t^2}}\right]\nonumber\\
    &=\frac{\tfrac{1}{2}\G\left(\tfrac{1}{4}\right)\G\left(\tfrac{1}{2}\right)}{\G\left(\tfrac{3}{4}\right)}\left(\frac{1}{4t^2}\right)^{3/4}\left[\l_4\left(\frac{1}{4t^2}\right)+e^{3\pi i/4}\l_4\left(1-\frac{1}{4t^2}\right)\right]\label{eq:numerator-i}.
  \end{align}
  Similarly, we obtain
  \begin{align}
    \F{\frac{1}{4},\frac{1}{4};\frac{1}{2};4t^2}    =\frac{\G\left(\tfrac{3}{4}\right)\G\left(\tfrac{1}{2}\right)}{\G\left(\tfrac{1}{4}\right)}\left(\frac{1}{4t^2}\right)^{1/4}\left[\l_4\left(\frac{1}{4t^2}\right)+e^{\pi i/4}\l_4\left(1-\frac{1}{4t^2}\right)\right]\label{eq:denominator-i}.
  \end{align}  
  Next, we divide \eqref{eq:numerator-i} by \eqref{eq:denominator-i} and use the formula for $\Omega_i$ to get
  \begin{align*}
    \frac{t}{2\pi\Omega_i^2} \frac{\F{\tfrac{3}{4},\tfrac{3}{4};\tfrac{3}{2};4t^2}}{\F{\tfrac{1}{4},\tfrac{1}{4};\tfrac{1}{2};4t^2}}
    &=\frac{\l_4\left(\frac{1}{4t^2}\right)+e^{3\pi i/4}\l_4\left(1-\frac{1}{4t^2}\right)}{\l_4\left(\frac{1}{4t^2}\right)+e^{\pi i/4}\l_4\left(1-\frac{1}{4t^2}\right)}.
  \end{align*}
  Therefore we have
  \[
    \frac{C_i(t)}{2\pi\Omega_i^2}=\frac{\l_4\left(\tfrac{1}{4t^2}\right)+e^{3\pi i/4}\l_4\left(1-\tfrac{1}{4t^2}\right)}{\l_4\left(\tfrac{1}{4t^2}\right)+ e^{\pi i/4}\l_4\left(1-\tfrac{1}{4t^2}\right)}.
  \]
  Taking $\g=1-1/4t^2$ and $\t$ as in Proposition \ref{thm:j-inv-q4} we get
  \begin{equation}\label{eq:tau-quotient-i}
    \frac{C_i(t)}{2\pi\Omega_i^2}=\frac{-i\sqrt{2}\t+e^{3\pi i/4}}{-i\sqrt{2}\t+e^{\pi i/4}}=\frac{(2\t-1)-i}{(2\t-1)-\ol{i}}=S_i(2\t-1).
  \end{equation}
  Therefore we obtain
  \[
    \frac{s_i^{-1}(C_i(t))+1}{2}=\t,
  \]
  and the inversion formula in Proposition \ref{thm:j-inv-q4} gives us
  \[
    j\left(\frac{s^{-1}_i(C_i(t))+1}{2}\right)=j(\t)=\frac{64\left(1+3\left(1-\tfrac{1}{4t^2}\right)\right)^3}{\left(1-\tfrac{1}{4t^2}\right)\left(\left(1-\tfrac{1}{4t^2}\right)-1\right)^2}=\frac{64(16t^2-3)^3}{4t^2-1}.
  \]

  For the case when $\arg(t)\not\in[0,\pi/2)$, then equation \eqref{eq:tau-quotient-i} has the form
  \[
    \frac{C_i(t)}{2\pi\Omega_i^2}=(-1)^a\frac{(2\t+(-1)^b)-i}{(2\t+(-1)^b)-\ol{i}}.
  \]
  The values of $a$ and $b$ as a function of the argument of $t$ are given by the following table:
  \[
    \begin{array}[!h]{|c|c|c|c|c|}\hline
      \arg(t) & [0,\pi/2) & [\pi/2,\pi) & [-\pi/2,0) & [-\pi,\pi/2)  \\ \hline
      a & 0 & 1 & 1 & 0\\ \hline
      b & 1 & 0 & 1 & 0\\ \hline
      \dfrac{C_i(t)}{2\pi\Omega_i^2} & \frac{(2\t-1)-i}{(2\t-1)-\ol{i}} & \frac{-\frac{1}{2\t+1}-i}{-\frac{1}{2\t+1}-\ol{i}} & \frac{-\frac{1}{2\t-1}-i}{-\frac{1}{2\t-1}-\ol{i}} & \frac{(2\t+1)-i}{(2\t+1)-\ol{i}}\\ [3mm] \hline 
      \frac{s_i^{-1}(C_i(t))+1}{2} & \t & \frac{\t-1}{2\t-1} & \frac{\t}{2\t+1} & \t+1 \\ [3mm] \hline
    \end{array}
  \]
  Clearly, all possible values of $\dfrac{s_i^{-1}(C_i(t))+1}{2}$ are $\mathrm{SL}_2(\Z)$-equivalent and thus their $j$-values are invariant.

\section{Proof of Theorem \ref{thm:sig6}}\label{sec:pf-thm1.2}

The following calculations depend on the argument of $t$. For the moment, we assume that $0<\arg(t)<\tfrac{\pi}{3}$. If we apply Proposition \ref{prop:linear-cov} to $\F{\tfrac{5}{6},\tfrac{5}{6};\tfrac{5}{3};-2t^{3}}$, using the functional equation $s\G(s)=\G(s+1)$ and the lambda notation in \eqref{eq:lambda} we get
  \begin{equation}
    \label{eq:numerator-i22}
    \F{\frac{5}{6},\frac{5}{6};\frac{5}{3};-2t^3}
    =\frac{\tfrac{2}{3}\G\left(\tfrac{1}{6}\right)\G\left(\tfrac{2}{3}\right)}{\G\left(\tfrac{5}{6}\right)}\left(\frac{1}{2t^3}\right)^{5/6}\Big[e^{5\pi i/6}\l_6(1/2t^{3})+\l_6(1+1/2t^{3})\Big]
  \end{equation}
  Similarly, if we apply Proposition \ref{prop:linear-cov} to $\F{\tfrac{1}{6},\tfrac{1}{6};\tfrac{1}{3};-2t^{3}}$ to get
  \begin{align}
    \F{\frac{1}{6},\frac{1}{6};\frac{1}{3};-2t^3}=\frac{\G\left(\tfrac{5}{6}\right)\G\left(\tfrac{1}{3}\right)}{\G\left(\tfrac{1}{6}\right)}\left(\frac{1}{2t^3}\right)^{1/6}\Big\lbrack e^{\pi i/6}\l_4(-1/2t^{3})+\l_6(1+1/2t^{3})\Big\rbrack\label{eq:denominator-i2}
  \end{align}
  We divide \eqref{eq:numerator-i22} by \eqref{eq:denominator-i2}, use Legendre's Duplication formula
    $\sqrt{\pi}\G(2s)=2^{2s-1}\G(s)\G\left(s+\tfrac{1}{2}\right),$
  and the definition of $\Omega_\r$, to obtain
  \[
    \frac{\cro(t)}{2\pi\Omega_\r^2}=\frac{e^{5\pi i/6}\l_6(-1/2t^{3})+\l_6(1+1/2t^{3})}{e^{\pi i/6}\l_6(-1/2t^{3})+\l_6(1+1/2t^{3})}
  \]
  Taking $\g=1/t$ and $\t$ as in Proposition \ref{thm:j-inv-q6} we get
  \begin{equation}\label{eq:tau-quotient-rho}
    \frac{\cro(t)}{2\pi\Omega_\r^2}=\frac{e^{5\pi i/6}-i\t}{e^{\pi i/6}-i\t}=\frac{\t-\r}{\t-\ol{\r}}=S_\r(\t).
  \end{equation}
  Therefore we obtain
  \[
    s_\r^{-1}(\cro(t))=\t,
  \]
  and the inversion formula in Proposition \ref{thm:j-inv-q6} gives
  \[
    j\left(s^{-1}_\r(\cro(t))\right)=j(\t)=-\frac{1728t^6}{(2t^3+1)}.
  \]

  For arbitrary arguments of $t$,  \eqref{eq:tau-quotient-rho} has the following three forms:
  \[
    \begin{array}[!h]{|c|c|c|c|}\hline
      \arg(t) & (-\tfrac{2\pi}{3},-\tfrac{\pi}{3})\cup(0,\pi/3)\cup(\tfrac{2\pi}{3},\pi)  & (-\pi,-\tfrac{2\pi}{3}]\cup(-\tfrac{\pi}{3},0]\cup(\tfrac{\pi}{3},\tfrac{2\pi}{3}]  & \{-\tfrac{\pi}{3},\tfrac{\pi}{3},\pi\}   \\ [2mm] \hline
      \dfrac{\cro(t)}{2\pi\Omega_\r^2} & \frac{\t-\r}{\t-\ol{\r}} & \frac{\t+1-\r}{\t+1-\ol{\r}} & \frac{-\frac{1}{\t-1}-\r}{-\frac{1}{\t-1}-\ol{\r}} \\ [5mm] \hline 
      s_\r^{-1}(\cro(t)) & \t & \t+1 & -\frac{1}{\t-1} \\ \hline
    \end{array}
  \]
  Recall that for the case $\arg(t)=\pi/3,\pi,-\pi/3$, we have to use Remark \ref{rm:identity-for-x>1}. Clearly, all possible values of $s_\r^{-1}(\cro(t)$ are $\mathrm{SL}_2(\Z)$-equivalent and thus their $j$-values are invariant.

\section*{Some comments}

Notice that in both the proofs above we did not require any conditions on $|t|$. However, we do require two conditions on $t$, namely $t$ must be in the domains where the inversion formulas in Propositions \ref{thm:j-inv-q4} and \ref{thm:j-inv-q6} are valid, and whenever $C_i(t)$ and $\cro(t)$ are well defined. The latter happens for $t\neq\pm1/2,$ and $t\neq(1/\sqrt[3]{2})\r^a$ for $a=1,3,5,$ respectively. The former happens when $\t=\t(t)\in\H$ which, by equations \eqref{eq:tau-quotient-i} and \eqref{eq:tau-quotient-rho}, is equivalent to $C_i(t)/2\pi\Omega_i^2\in\mathbb{D}$ and $\cro(t)/2\pi\Omega_\r^2\in\mathbb{D}$ for Theorems \ref{thm:sig4} and \ref{thm:sig6} respectively. This can be verified easily using the Maximum Modulus Principle on $|t|<1/2$ and $|t|<1/\sqrt[3]{2}$ respectively.

Furthermore, from the definition of $_2F_1$, $C_i(t)$ and $\cro(t)$ are discontinuous on the rays
\[
    \{u(-1)^a\mid u\geq 1/2, a=0,1\}\quad\text{and}\quad \{u\r^a\mid u\geq 1/\sqrt[3]{2}, a=1,3,5\}
\]
respectively, but, as apparent from the Tables above $j((s_i^{-1}C_i(t)+1)/2)$ and $j(s_\r^{-1}\cro(t))$ become continuous (and thus analytic) functions of $t$ which can be proved  using Remark \ref{rm:identity-for-x>1}. Therefore, Theorems \ref{thm:sig4} and \ref{thm:sig6} are valid for all $t$ except for two and three points respectively.

\section{Examples}\label{sec:examples}
Here we offer some examples.
\begin{example}
    It is well known that $j(\sqrt{-2})=8000$. We verify this using Theorem \ref{thm:sig4}. First we solve the degree six equation in $t$
    \[
        64\frac{(16t^2-3)^3}{4t^2-1}=8000.
    \]
    One solution is
    \[
        t_0=\frac{i}{4\sqrt{2}}\sqrt{5\sqrt{2}-1}=i\cdot 0.4355695915...,
    \]
    so that $|t_0|<1/2$. We approximate $C_i(t_0)$ using the first 3000 terms of its power series expansion to get:
    \[
        C_i(t_0)=i\cdot 0.375476877103748...
    \]
    Thus
    \[
        \t_0:=\frac{s_i^{-1}(C_i(t_0))+1}{2}=0.333333333333333... + i\cdot 0.471404520791031...
    \]
    Notice that $\t_0\neq\sqrt{-2}$, but they are $\SL_2(\Z)$-equivalent. Indeed, we have
    \[
        -\frac{1}{\t_0}+1=i\cdot 1.414213562373095...\approx\sqrt{-2}.
    \]
    In fact, the above approximation is correct up to 364 decimal places.
\end{example}
\begin{example}
    
    Now we use Theorem \ref{thm:sig6} to verify $j((1+\sqrt{-7})/2)=-3375$. We solve the degree 6 equation
    \[
        -\frac{1728 t^6}{2t^3+1}=-3375.
    \]
    The solutions satisfy
    \[
        t^3=\frac{5}{64}(25\pm 3\sqrt{105}),
    \]
    so we take the real cubic root of $\tfrac{5}{64}(25-3\sqrt{105})$, which is
    \[
        t_0=\sqrt[3]{\frac{5}{64}(25-3\sqrt{105})}=-0.765459354046599...
    \]
    Since $|t_0|<1/\sqrt[3]{2}\approx 0.793700...$, we can approximate $\cro(t_0)$ using the first 3000 terms of the power series expansion :
    \[
        \cro(t_0)=0.538697866211295...
    \]
    Thus
    \[
        \t_0:=s_\r^{-1}(\cro(t_0))=0.500000000000000... + i\cdot 1.322875655532295...\approx\frac{1+\sqrt{-7}}{2}.
    \]
    In fact, the above approximation is correct up to 145 decimal places.
\end{example}

\begin{example}
    To illustrate the general inversion process, we try to find $\t_0$ such that $j(\t_0)=-50,000.$ We solve the degree 6 equation
    \[
        -\frac{1728 t^6}{2t^3+1}=-50,000.
    \]
    The solutions satisfy
    \[
        t^3=\frac{25}{108}(125-\sqrt{16165}),
    \]
    so we take the real cubic root of $\tfrac{25}{108}(125-\sqrt{16165})$, which is
    \[
        t_0=\sqrt[3]{\frac{25}{108}(125-\sqrt{16165})}=-0.791446942386710....
    \]
    Notice that $|t_0|<1/\sqrt[3]{2}\approx 0.793700...$, but it is very close to the upper bound, which means that we have to use more terms in the power series of $\cro(t)$ to get a reasonable approximation. Using the first 5000 terms yields
    \[
        \cro(t_0)=0.855243324301038...
    \]
    and
    \[
        \t_0=0.500000000000000...+i\cdot 1.724359831532281....
    \]
    We find that
    \[
        j(\t_0)=-49,999.9999999999999996....
    \]
    In fact, the above approximation is correct up to 16 decimal places when we just use the first 5000 terms of power series.
\end{example}

\bibliography{bibliography}
\bibliographystyle{alpha}

\end{document}